\title
{Flows on Bidirected Graphs}
\author{
    Matt DeVos\\
    {\normalsize mdevos@sfu.ca}
}

\date{}

\documentclass[12pt]{article}

\usepackage{graphicx}
\usepackage{amsfonts}
\usepackage{amssymb}
\usepackage{amsmath}
\usepackage{amsthm}
\usepackage{latexsym}

\begin{document}

\bibliographystyle{plain}
\maketitle
\setcounter{page}{1}
\newtheorem{theorem}{Theorem}[section]
\newtheorem{lemma}[theorem]{Lemma}
\newtheorem{corollary}[theorem]{Corollary}
\newtheorem{proposition}[theorem]{Proposition}
\newtheorem{definition}[theorem]{Definition}
\newtheorem{claim}{Claim}
\newtheorem{conjecture}[theorem]{Conjecture}
\newtheorem{observation}[theorem]{Observation}

\begin{abstract}
The study of nowhere-zero flows began with a key observation of Tutte
that in planar graphs, nowhere-zero k-flows are dual to k-colourings (in
the form of k-tensions).  Tutte conjectured that every graph without a
cut-edge has a nowhere-zero 5-flow.  Seymour proved that every such graph
has a nowhere-zero 6-flow.

For a graph embedded in an orientable surface of higher genus, flows are
not dual to colourings, but to local-tensions.  By Seymour's theorem, every
graph on an orientable surface without the obvious obstruction has a
nowhere-zero 6-local-tension.
Bouchet conjectured that the same should hold true on non-orientable
surfaces.  Equivalently, Bouchet conjectured that every bidirected graph
with a nowhere-zero $\mathbb{Z}$-flow has a nowhere-zero 6-flow.  
Our main result establishes that every such graph has a nowhere-zero 12-flow.
\end{abstract}

\section{Introduction}

Throughout the paper, we consider only finite graphs, which may have
loops and parallel edges.  Let $G$ be a graph and let $X \subseteq V(G)$.  We let
$\delta_G(X)$ denote the set of edges with exactly one endpoint in $X$, and call any 
such set an \emph{edge cut}.  When $X = \{x\}$ we abbreviate this notation to $\delta_G(x)$ 
and we drop the subscript $G$ when the graph is clear from context.  

A {\it signature} of $G$ is a function $\sigma : E(G) \rightarrow \{ \pm 1 \}$.  We say that an edge $e \in E(G)$
is {\it positive} if $\sigma(e) = 1$ and {\it negative} if $\sigma(e) = -1$.  
For $S \subseteq E(G)$, we let $\sigma(S) = \prod_{e \in S} \sigma(e)$ 
and for a subgraph $H \subseteq G$, we let $\sigma(H) = \sigma(E(H))$.  A cycle $C \subseteq G$
is  \emph{balanced (with respect to $\sigma$)} if $\sigma(C) = 1$ and \emph{unbalanced (with respect to $\sigma$)} 
if $\sigma(C) = -1$.  

Let $v \in V(G)$, and modify $\sigma$ to make a new signature
$\sigma'$ by changing $\sigma'(e) = - \sigma(e)$ for every $e \in
\delta(v)$.  We say that $\sigma'$ is obtained from $\sigma$ by
making a {\it flip} at the vertex $v$ and we define two signatures of $G$ 
to be \emph{equivalent} if one can be obtained from the other by a sequence
of flips.  It is a straightforward exercise to verify that the following statements 
are equivalent for any two signatures $\sigma_1, \sigma_2$ of $G$.
\begin{itemize}
\item $\sigma_1$ and $\sigma_2$ are equivalent.
\item There is an edge cut $S$ so that $\sigma_1$ and $\sigma_2$ differ precisely on $S$.
\item Every cycle $C \subseteq G$ satisfies: $C$ is balanced with respect to $\sigma_1$ if and only if it 
	is balanced with respect to $\sigma_2$.  
\end{itemize}
A {\it signed graph} consists of a graph $G$ equipped with a signature $\sigma_G$.  We say that $G$
is \emph{balanced} if every cycle of $G$ is balanced (or equivalently its $\sigma_G$ is equivalent
to the constant 1).  

Following Bouchet \cite{Bo}, we will treat each edge of the graph $G$ as composed of two half edges.  
So, every half edge $h$ is contained in exactly one edge, denoted $e_h$, and incident with 
exactly one vertex which must be an end of $e_h$.  A non-loop edge contains one half edge incident 
with each end, while a loop edge contains two half-edges each incident with the only end.  We let $H(G)$ 
denote the set of half edges in $G$.  For every $v \in V(G)$ we let $H(v)$ denote the set of half edges 
incident with $v$, and for every $e \in E(G)$ we let $H(e)$ denote the set of half edges contained in $e$.  

\begin{figure}[h]
  \centering
    \includegraphics[height=2.4cm]{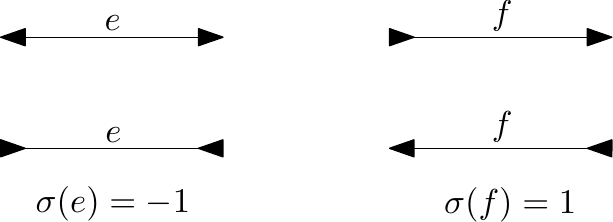}
  \label{fig1}
  \caption{Orientations of signed edges}
\end{figure}

An \emph{orientation} of a signed graph $G$ is a function $\tau : H(G) \rightarrow \{ \pm 1 \}$ with the property 
that $\prod_{h \in H(e)} \tau(h) = -\sigma_G(e)$ for every $e \in E(G)$.  If $h$ is a half edge incident with the 
vertex $v$, then $h$ is directed \emph{toward} $v$ if $\tau(h) = 1$ and directed \emph{away} if $\tau(h) = -1$ (see Figure \ref{fig1}).  
A \emph{bidirected} graph consists of a signed graph $G$ together with an orientation $\tau_G$.  When the signature is the 
constant 1, each edge contains two half edges which are consistently oriented, and this aligns with the usual 
notion of an orientation of an (ordinary) graph.  Indeed we will view ordinary digraphs and bidirected graphs with all edges
positive as the same.  


Let $G$ be a bidirected graph, let $\Gamma$ be an abelian group,
and let $\phi : E(G) \rightarrow \Gamma$ be a function.  We define the {\it boundary}
of $\phi$ be the function $\partial \phi : V(G) \rightarrow \Gamma$ given by the rule
\[ \partial \phi (v) = \sum_{h \in H(v)} \tau_G(h) \phi(e_h). \]
We define $\phi$ to be a {\it flow} if $\partial \phi = 0$; if in addition $\Gamma = \mathbb{Z}$ and
$|\phi(e)| < k$ for every $e \in E(G)$ we call $\phi$ a $k$-\emph{flow}.  If $0 \not\in \phi(E(G))$ we will say that 
$\phi$ is {\it nowhere-zero}.

Suppose that $\phi$ is a flow of a bidirected graph $G$ with orientation $\tau$ and let $e \in E(G)$.  Now, modify $\tau$ to form a new
orientation $\tau'$ by changing $\tau'(h) = -\tau(h)$ for every $h \in H(e)$ and modify $\phi$ to form a new function 
$\phi'$ by changing $\phi'(e) = -\phi(e)$.  After these adjustments, $\phi'$ is a flow of the new bidirected graph 
given by $G$ and $\tau'$.  Furthermore, $\phi'$ is nowhere-zero if and only if $\phi$ is nowhere-zero and $\phi'$ is a k-flow if and
only if $\phi$ is a k-flow.  So, as in the case of ordinary graphs, the existence of a nowhere-zero $\Gamma$ flow 
or nowhere-zero $k$-flow in a bidirected graph does not depend on the orientation.  Accordingly, we say that a signed graph 
has a \emph{nowehere-zero} $\Gamma$-flow or $k$-flow if some (and thus every) orientation of it has such a flow.  

Again consider a flow $\phi$ in a bidirected graph given by $G$ together with $\sigma$ and $\tau$, and now let 
$v \in V(G)$.    Modify $\sigma$ to make a new signature $\sigma'$ by making a flip at the vertex $v$, and
modify $\tau$ to form $\tau'$ by changing $\tau'(h) = - \tau(h)$ for every $h \in H(v)$.  It follows that $\phi$ is a flow 
of the bidirected graph given by $G$ together with $\sigma'$ and $\tau'$.  Therefore, the existence of a nowhere-zero 
$\Gamma$-flow or $k$-flow in a signed graph is invariant under changing the signature to an equivalent one.  

The following conjecture of Bouchet is the central question in the study of nowhere-zero flows in signed graphs.

\begin{conjecture}[Bouchet's 6-Flow Conjecture \cite{Bo}]
\label{6FC}
Every signed graph with a nowhere-zero $\mathbb{Z}$-flow has a nowhere-zero $6$-flow.
\end{conjecture}

Bouchet \cite{Bo} proved that the above conjecture holds with 6 replaced by 216, and gave an example to show
that 6, if true, would be best possible.  Zyka and independently Fouquet proved that the above conjecture is true with 
6 replaced by 30.  Our main result is as follows.

\begin{theorem}
\label{12FT}
Every signed graph with a nowhere-zero $\mathbb{Z}$-flow has a
nowhere-zero $12$-flow.
\end{theorem}

For 4-edge-connected signed graphs, Bouchet's conjecture holds true.  In fact, every 4-edge-connected signed graph 
with a nowhere-zero $\mathbb{Z}$-flow also has a nowhere-zero 4-flow.  This was proved recently by Raspaud and 
Zhu \cite{RZ}.  It was also proved earlier by this author (using essentially the same argument) in his thesis \cite{De}, but never published.

Since we will be focused on graphs with nowhere-zero $\mathbb{Z}$-flows, it will be helpful to have a precise description
of this class of graphs.  The following characterization follows from Bouchet's original paper, although it is not explicitly 
stated there.

\begin{proposition}[Bouchet]
\label{nzzflow}
A connected signed graph $G$ does not have a nowhere-zero $\mathbb{Z}$-flow if and only if one of the following holds:
\begin{enumerate}
\item $\sigma_G$ is equivalent to a signature with exactly one negative edge.
\item $G$ has a cut-edge $e$ for which $G \setminus e$ has a balanced component.
\end{enumerate}
\end{proposition}

\section{Modular Flows on Bidirected Graphs}

In this section we reduce the problem of finding a nowhere-zero $12$-flow to that of finding a certain
type of flow in the group $\mathbb{Z}_2 \times \mathbb{Z}_3$.  
If $\phi : S \rightarrow X_1 \times X_2 \times \ldots \times X_n$ we will
let $\phi_i$ denote the projection of $\phi$ onto $X_i$.
If $G$ is a bidirected graph and $\phi : E(G) \rightarrow \mathbb{Z}_2 \times \mathbb{Z}_3$ is a flow,
we will say that $\phi$ is {\it balanced} if $\sigma_G( {\mathit supp}(\phi_1) ) = 1$.  
The purpose of this section is to prove the following lemma.

\begin{lemma}
\label{bnzl}
Every signed graph with a balanced nowhere-zero $\mathbb{Z}_2 \times \mathbb{Z}_3$-flow has a 
nowhere-zero 12-flow.
\end{lemma}

Our proof of Lemma \ref{bnzl} is based on a theorem of Bouchet concerning chain groups which 
simplifies to the following in the case of bidirected graphs.

\begin{theorem}[Bouchet \cite{Bo}]
\label{T42}
Let $G$ be a bidirected graph, let $\phi$ be a $\mathbb{Z}$-flow of $G$, and
let $k > 0$.  Then there exists a 2k-flow $\phi'$ of $G$ so that
$\phi'(e) \equiv \phi(e)$ (mod $k$) for every $e \in E(G)$.
\end{theorem}

We proceed toward the proof of Lemma \ref{bnzl} with some straightforward lemmas.
We will repeatedly call upon the following simple formula which holds for every bidirected graph $G$ 
and function $\phi : E(G) \rightarrow \Gamma$.
\begin{equation}
\label{sumboundary}
\sum_{v \in V(G)} \partial \phi (v) = \sum_{e \in E(G) : \sigma_G(e) = -1} \pm 2 \phi(e)
\end{equation}

\begin{lemma}
\label{posadj}
Let $G$ be a connected bidirected graph with all edges positive and let $\Gamma$ be an abelian group.  
If $\mu : V(G) \rightarrow \Gamma$ satisfies $\sum_{v \in V(G)} \mu(v) = 0$, there exists $\phi : E(G) \rightarrow \Gamma$ with 
$\partial \phi = \mu$.
\end{lemma}

\noindent{\bf Proof:} Let $T$ be a spanning tree of $G$.  We begin with $\phi$ identically $0$ and modify it in steps as follows.  
Choose  a leaf vertex $v$ of $T$ incident with the leaf edge $e \in E(T)$, adjust $\phi(e)$ so that $\partial \phi(v) = \mu(v)$, and
then modify $T$ by delting $v$.  When $T$ consists of a single vertex $u$, the function $\phi$ satisfies $\partial \phi(v) = \mu(v)$ for
every $v \in V(G) \setminus \{u\}$.  Our assumptions and equation (\ref{sumboundary}) imply that
$\sum_{v \in V(G)} \mu(v) = 0 = \sum_{v \in V(G)} \partial \phi(v)$ so $\partial \phi = \mu$ as desired.
\quad\quad$\Box$

\begin{lemma}
\label{L43}
Let $G$ be a connected bidirected graph with an unbalanced
cycle, and let $\mu : V(G) \rightarrow \mathbb{Z}$ be a function with
$\sum_{v \in V(G)} \mu(v)$ even.  Then there exists a function
$\phi : E(G) \rightarrow \mathbb{Z}$ such that $\partial \phi = \mu$.
\end{lemma}

\noindent{\bf Proof:}
Let $C$ be an unbalanced cycle of $G$, and let $u \in V(C)$ and $e \in E(C)$
be incident.  Choose a spanning tree $T \subseteq G$ so that
$C \setminus e \subseteq T$.  As in the preceding lemma we may  
choose a function $\phi : E(G) \rightarrow \mathbb{Z}$ supported on a subset of $E(T)$ so that 
$\partial \phi(v) = \mu(v)$ for every $v \in V(T) \setminus \{u\}$.  It follows from our
assumptions and equation (\ref{sumboundary}) that $\mu(u) - \partial \phi(u)$ is even.  
Since $C$ is unbalanced, we may choose a function $\psi : E(G) \rightarrow \{ -1,0,1\}$ with support 
$E(C)$ so that
\[ \partial \psi(v) = \left\{ \begin{array}{ll}
    2   &   \mbox{if $v=u$} \\
    0   &   \mbox{otherwise}
    \end{array} \right. \]
Now the function $\phi + \frac{\mu(u) - \partial \phi(u)}{2} \psi$ has boundary $\mu$ as desired.
\quad\quad$\Box$

\begin{lemma}
\label{L44}
Let $G$ be a connected bidirected graph, let $p$ be a prime, let $\psi$ be
a $\mathbb{Z}_p$-flow of $G$, and assume that either $p$ is odd or that
$\sigma_G({\mathit supp}(\psi)) = 1$.  Then there is a $\mathbb{Z}$-flow
$\phi$ of $G$ so that $\phi(e) \equiv \psi(e)$ $(mod \;\, p)$ for every
$e \in E(G)$.
\end{lemma}

\noindent{\bf Proof:}
Choose $\phi : E(G) \rightarrow \mathbb{Z}$ so that
$\phi(e) \equiv \psi(e)$ (mod $p$) for every $e \in E(G)$.  Since
$\psi$ is a $\mathbb{Z}_p$-flow, we will have
$\partial \phi (v)$ a multiple of $p$ for every $v \in V(G)$.  It follows from equation (\ref{sumboundary}) that 
$\sum_{v \in V(G)} \partial \phi (v)$ is even.  Furthermore, if $p=2$ then we have an even number of negative edges (by assumption),
so in this case $\sum_{v \in V(G)} \partial \phi (v)$ is a multiple of $4$.  By the preceding lemma, we may choose
$\eta : E(G) \rightarrow \mathbb{Z}$ so that $\partial \eta = (1/p) \partial \phi$.  Now $\phi - p \eta$ is a flow with the desired properties.
\quad\quad$\Box$

\bigskip

We are now ready to prove Lemma \ref{bnzl}

\bigskip

\noindent{\bf Proof of Lemma \ref{bnzl}}
Let $G$ be a bidirected graph with a balanced nowhere-zero
$\mathbb{Z}_2 \times \mathbb{Z}_3$-flow $\psi$.
By Lemma \ref{L44} we may choose integer flows $\phi_1,\phi_2$ so that
$\phi_1 \equiv \psi_1$ (mod 2) and $\phi_2 \equiv \psi_2$
(mod 3).  Now $\eta = 3 \phi_1 + 2 \phi_2$ is an integer flow with 
the property that $\eta(e) \not\equiv 0$
(mod 6) for every $e \in E(G)$.  By Theorem \ref{T42} we may now choose
an integer flow $\eta'$ so that $\eta'(e) \equiv \eta(e)$ (mod 6)
and $|\eta'(e)| < 12$ for every $e \in E(G)$.  Now $\eta'$ is a
nowhere-zero 12-flow as desired.
\quad\quad$\Box$

\section{Seymour's 6-Flow Theorem}

Thanks to the previous section, our remaining task is to show that every signed graph with a nowhere-zero
$\mathbb{Z}$-flow also has a balanced nowhere-zero $\mathbb{Z}_2 \times \mathbb{Z}_3$-flow.  Our proof of this will follow one
of Seymour's proofs of the $6$-Flow Theorem.  However, our notation and process are 
sufficiently complicated to obscure this link.  To make this connection apparent and to assist the 
reader in understanding our forthcoming argument, we next give a sketch of Seymour's proof
using similar language to that appearing in ours.

\begin{theorem}[Seymour]
Every 2-edge-connected graph has a nowhere-zero $\mathbb{Z}_2 \times \mathbb{Z}_3$-flow.
\end{theorem}

\noindent{\bf Proof Sketch:} 
The first step is a straightforward reduction to 3-edge-connected graphs.  Suppose that $\{e,f\}$ is a 2-edge cut of the oriented graph $G$ and 
consider $G/e$.  If $\phi$ is a nowhere-zero $\mathbb{Z}_2 \times \mathbb{Z}_3$-flow of $G/e$ then it follows from basic principles that $\phi$ may be extended to a flow of $G$ by assigning some value to the edge $e$.  However, since $\{e,f\}$ is an edge cut 
in $G$ we must have $\phi(e) = \pm \phi(f)$ and thus $\phi$ is nowhere-zero.  So to prove the theorem, it suffices to consider 3-edge-connected graphs.

Next we reduce to cubic graphs.  If $G$ is an oriented 3-edge-connected graph with a vertex $v$ of degree at 
least 4, then we may uncontract an edge at $v$ so that the resulting graph $G'$ remains 3-edge-connected
(this is a straightforward exercise).  Now if $\phi'$ is a nowhere-zero $\mathbb{Z}_2 \times \mathbb{Z}_3$-flow 
of $G'$ then $\phi = \phi'|_{E(G)}$ is a nowhere-zero $\mathbb{Z}_2 \times \mathbb{Z}_3$-flow of $G$.  
So, to prove Seymour's theorem, it suffices to consider 3-connected cubic graphs.  

Now we will prepare for an inductive argument.  First note that if $H$ is a proper induced subgraph
of a 3-connected cubic graph, then $H$ is a subcubic graph and for every $X \subseteq V(H)$ we have
\[ (\star) \quad\quad  |\delta_H(X)| + \sum_{x \in X} (3- \mathit{deg}_H(x)) \ge 3.  \]
For an oriented subcubic graph $H$ let us define a \emph{watering} to be a function 
$\phi : E(H) \rightarrow \mathbb{Z}_2 \times \mathbb{Z}_3$ with the following property:
\[ \partial \phi (v) = \left\{ \begin{array}{ll}
    (0,0)       &   \mbox{if ${\mathit deg}(v) = 3$}    \\
    (0,\pm1)    &   \mbox{if ${\mathit deg}(v) = 1,2$}
    \end{array} \right. \]
This more general concept will permit us to work inductively on subgraphs.  
Let $G$ be an oriented 3-connected cubic graph, let $u \in V(G)$, and let $\delta(u) = \{e_1,e_2,e_3\}$.  Suppose that $\phi$ is a nowhere-zero
watering of $G \setminus u$.  Then we may extend $\phi$ to have domain $E(G)$ by choosing $\phi(e_i) = (0,\pm 1)$ in such a way 
that $\partial \phi(v) = 0$ for all $v \in V(G) \setminus \{u\}$.  It now follows that $\phi$ is a nowhere-zero flow of $G$.  So, to complete 
the proof, it suffices to show that every subcubic graph satisfying $(\star)$ has a nowhere-zero watering.  

We shall prove that every subcubic graph $H$ satisfying $(\star)$ has a nowhere-zero watering by induction on 
$|V(H)|$.  If $H$ has a vertex of degree $1$ then by induction $H \setminus v$ has a nowhere-zero watering, and this 
may be extended to a nowhere-zero watering of $H$.  So, we may assume that $H$ has minimum degree $2$.  Clearly 
we may also assume $H$ is connected.

We claim that $H$ must have a cycle containing at least two vertices of degree 2, and let us call any such cycle \emph{removable}.  
If $H$ is 2-connected, there exist at least two vertices of degree 2 (by applying $(\star)$ to $V(H)$), and thus we have a 
removable cycle.  Otherwise we may choose  leaf-block of $H$ and by $(\star)$ this leaf block must have at least two degree 2 
vertices, so again we find a removable cycle.  

Now choose a removable cycle $C$ and consider the graph $H' = H \setminus V(C)$.  This is another subcubic graph
satisfying $(\star)$, so by induction we may choose a nowhere-zero watering $\phi'$ of $H'$.  We extend $\phi'$ to 
a function $\phi : E(H) \rightarrow \mathbb{Z}_2 \times \mathbb{Z}_3$ as follows.
\[ \phi(e) = \left\{ \begin{array}{ll}
    (0,\pm 1)&	\mbox{if $e \in \delta(V(C))$} \\
    (1,0)       &   \mbox{if $e \in E(C)$}          \\
    (0,1)       &   \mbox{if $e$ is a chord of $C$} 
    \end{array} \right. \]
Since every vertex in $H'$ incident with an edge in $\delta(C)$ has degree at most two in $H'$ we may choose the values $\phi(e)$
on edges $e \in \delta(V(C))$ so that $\phi$ satisfies the boundary condition for a watering at every vertex in $H'$.  Note that by construction
the function $\phi$ satisfies $\partial \phi_1 (v) = 0$ for every $v \in V(C)$.  So, we only need adjust $\partial \phi_2(v)$ for vertices $v \in V(C)$ 
to obtain a watering.  For every $v \in V(C) \cap V_2(H)$ let $\beta_v$ be a variable in $\mathbb{Z}_3$.  Since $C$ 
was removable, there are at least two such vertices, so we may choose $\pm 1$ assignments to the $\beta_v$ variable so that the following equation is satisfied
\[ \sum_{v \in V(C)} \partial \phi_2(v) = \sum_{v \in V(C) \cap V_2(H)} \beta_v. \]
By Lemma \ref{posadj} we may choose a function $\psi : E(C) \rightarrow \mathbb{Z}_3$ so that 
\[ \partial \psi(v) = \left\{ \begin{array}{cl}
					\beta_v 	- \partial \phi_2(v) 		&	\mbox{if $v \in V(C) \cap V_2(G)$}	\\
							- \partial \phi_2(v) 		&	\mbox{if $v \in V(C) \setminus V_2(G)$}	
					\end{array} \right. \]		
Now modify $\phi$ by adding $\psi$ to $\phi_2$.  The resulting function is a nowhere-zero watering of $H$ as desired.
\quad\quad$\Box$

\section{Restricted Flows in Digraphs}

In this section we will prove an elementary result concerning flows in 
ordinary graphs.  This will be used to reduce the problem of finding 
a balanced nowhere-zero $\mathbb{Z}_2 \times \mathbb{Z}_3$-flow to graphs with some
better connectivity properties.

If $G$ is a directed graph and $X \subseteq V(G)$, we let $\delta^+(X)$ denote the set of edges with initial vertex in $X$ and
terminal vertex in $V(G) \setminus X$.  We let $\delta^-(X) = \delta^+( V(G) \setminus X )$.
The goal of this section is to prove the following lemma, which will be used to build up the connectivity required for the proof of our 12-flow theorem.

\begin{lemma}
\label{L31}
Let $G$ be a directed graph, let $\Gamma$ be an abelian group, and assume
that $G$ has a nowhere-zero $\Gamma$-flow.  If $u \in V(G)$ is a vertex with
${\mathit deg}(u) \le 3$ and
$\gamma : \delta(u) \rightarrow \Gamma \setminus \{0\}$ satisfies
$\partial \gamma (u) = 0$, then there is a nowhere-zero $\Gamma$-flow
$\phi$ of $G$ so that $\phi|_{\delta(u)} = \gamma$.
\end{lemma}

After a few definitions, we will prove a lemma of Seymour, from which the
above lemma will easily follow.  Let $G$ be a directed graph, let $T \subseteq E(G)$, and let $\Gamma$ be an
abelian group.  For any function
$\gamma : T \rightarrow \Gamma$,
we will let ${\mathcal F}_{\gamma}(G)$ denote the number of nowhere-zero
$\Gamma$-flows $\phi$ of $G$ with $\phi(e) = \gamma(e)$ for every $e \in T$.
For every $X \subseteq V(G)$, let $\alpha_X : E(G) \rightarrow \{-1,0,1\}$ be
given by the rule
\[ \alpha_X(e) = \left\{ \begin{array}{cl}
    +1  &   \mbox{if $e \in \delta^+(X)$}   \\
    -1  &   \mbox{if $e \in \delta^-(X)$}   \\
    0   &   \mbox{otherwise}
    \end{array} \right. \]
If $\gamma_1,\gamma_2 : T \rightarrow \Gamma$, we will call
$\gamma_1,\gamma_2$
{\it similar} if for every $X \subseteq V(G)$, it holds that
\begin{eqnarray}
\sum_{e \in T} \alpha_X(e) \gamma_1(e) = 0
    & \mbox{if and only if} &
\sum_{e \in T} \alpha_X(e) \gamma_2(e) = 0
\end{eqnarray}

\begin{lemma}[Seymour - personal communication]
\label{L32}
Let $G$ be a directed graph and let $T \subseteq E(G)$.
If $\gamma_1,\gamma_2 : T \rightarrow \Gamma$
are similar, then ${\mathcal F}_{\gamma_1}(G) = {\mathcal F}_{\gamma_2}(G)$.
\end{lemma}

\noindent{\bf Proof:}
We proceed by induction on the number of edges in $E(G) \setminus T$.
If this set is empty, then
${\mathcal F}_{\gamma_i}(G) \le 1$ and ${\mathcal F}_{\gamma_i}(G) = 1$
if and only if $\gamma_i$ is a flow of $G$ for $i = 1,2$.  Thus, the result
follows by the assumption.  Otherwise, choose an edge $e \in E(G) \setminus T$.
If $e$ is a cut-edge then ${\mathcal F}_{\gamma_i}(G) = 0$ for $i=1,2$.  If
$e$ is a loop, then we have inductively that
\[ {\mathcal F}_{\gamma_1}(G) =
(|\Gamma| -1 ) {\mathcal F}_{\gamma_1}(G \setminus e) =
(|\Gamma| -1 ) {\mathcal F}_{\gamma_2}(G \setminus e) =
{\mathcal F}_{\gamma_1}(G)  \]
Otherwise, applying induction to $G \setminus e$ and $G / e$ we have
\[ {\mathcal F}_{\gamma_1}(G) = {\mathcal F}_{\gamma_1}(G / e)
- {\mathcal F}_{\gamma_1}(G \setminus e)
= {\mathcal F}_{\gamma_2}(G / e) - {\mathcal F}_{\gamma_2}(G \setminus e) =
{\mathcal F}_{\gamma_2}(G).  \quad\quad\Box \]

\noindent{\bf Proof of Lemma \ref{L31}}
If $\phi$ is a nowhere-zero $\Gamma$-flow of $G$, then
$\phi|_{\delta(v)}$ is similar to $\gamma$.  Thus by Lemma \ref{L32}, we have
that ${\mathcal F}_{\gamma}(G) = {\mathcal F}_{\phi|_{\delta(v)}}(G) \neq 0$.
\quad\quad$\Box$

\section{Reductions}

Just as Seymour's 6-flow theorem reduced to a problem on a certain type of subcubic graph, so shall our problem of finding a balanced nowhere-zero 
${\mathbb Z}_2 \times {\mathbb Z}_3$-flow.  Next we introduce this family of graphs and then give this reduction.  

We define a {\it shrubbery} to be a signed graph $G$ with the following properties:
\begin{enumerate}
\item $G$ has maximum degree at most $3$.
\item If $H \subseteq G$ is a component of $G$ and every vertex in $H$ has degree three, then $H \setminus e$ contains an unbalanced cycle for every $e \in E(H)$. 
\item For every $X \subseteq V(G)$ with $|X| \ge 2$, if $G[X]$ is balanced then \[ | \delta(X) | + \sum_{x \in X} (3 - {\mathit deg}(x))  > 3 \]
\item $G$ has no balanced cycles of length 4
\end{enumerate}
Fix an orientation of $G$, and define a {\it watering} to be a function $\phi: E(G) \rightarrow \mathbb{Z}_2 \times \mathbb{Z}_3$ so that
\[ \partial \phi (v) = \left\{ \begin{array}{ll}
    (0,0)       &   \mbox{if ${\mathit deg}(v) = 3$}    \\
    (0,\pm 1 )  &   \mbox{if ${\mathit deg}(v) = 1,2$}
    \end{array} \right. \]
Note that as in the case of flows, the chosen orientation does not effect the existence of a nowhere-zero 
watering.  Similarly, these properties are unaffected when we replace $\sigma_G$ by an equivalent signature. 
Our main lemma on shrubberies appears next.  We define a \emph{theta} to be a graph consisting of two vertices and three internally 
disjoint paths between them.  A theta is \emph{unbalanced} if it is not balanced (i.e it contains an unbalanced cycle).

\begin{lemma}
\label{workhorse}
Every shrubbery has a nowhere-zero watering.  Furthermore, if $G$ is a shrubbery with an unbalanced theta or loop and $\epsilon = \pm 1$, 
then $G$ has a nowhere-zero watering $\phi$ for which $\sigma_G ( {\mathit supp}( \phi_1 ) ) = \epsilon$.
\end{lemma}

Lemma \ref{workhorse} will be proved in the final section of the paper.  In the remainder of this section, we will 
show that it implies Theorem \ref{12FT}.  This requires the following.

\begin{lemma}
\label{lred}
Let $G$ be a signed graph with a nowhere-zero $\mathbb{Z}$-flow.  
\begin{enumerate}
\item If $e \in E(G)$ is positive, then $G/e$ has a nowhere-zero $\mathbb{Z}$-flow.
\item If $v \in V(G)$ has $\mathit{deg}(v) \ge 4$, we may uncontract a positive edge $e$ at $v$ forming  
vertices $v_1,v_2$ so that $\mathit{deg}(v_i) \ge 3$ for $i=1,2$ and the new graph has a nowhere-zero 
$\mathbb{Z}$-flow.
\end{enumerate}
\end{lemma}

\noindent{\bf Proof:} Fix an orientation $\tau$ of $G$ and let $\phi : E(G) \rightarrow \mathbb{Z}$ be a 
nowhere-zero flow.  The first part follows from the observation that $\phi|_{E(G) \setminus \{e\}}$ is a 
nowhere-zero flow of $G/e$.  For the second part, consider the equation 
$\sum_{h \in H(v)} \tau(h) \phi( e_h ) = 0$.  Since there are at least four terms in this sum
and all are nonzero, we may partition $H(v)$ into $\{H_1,H_2\}$ so that $|H_i| \ge 2$ 
and $\sum_{h \in H_i} \tau(h) \phi(e_h)  \neq 0$ for $i=1,2$.  Now form a new bidirected graph $G'$ from $G$ 
by uncontracting a positive edge $e$ at $v$ to form two new vertices $v_1$ and $v_2$ where $v_i$ is incident
with the half edges in $H_i$ and one from $e$.  Giving $e$ an arbitrary orientation, we may then extend $\phi$ 
to a nowhere-zero flow in $G'$ as desired.
\quad\quad$\Box$

\bigskip

We are now ready to prove Theorem \ref{12FT} using Lemma \ref{workhorse}.

\bigskip

\noindent{\bf Proof of Theorem \ref{12FT}:} We shall prove that every signed graph $G$ with a nowhere-zero $\mathbb{Z}$-flow 
has a nowhere-zero $12$-flow by induction on $\sum_{v \in V(G)} | {\mathit deg}(v) - 5/2 |$.  By the inductive hypothesis, we may 
assume that $G$ is connected.  By Seymour's 6-Flow Theorem, we may also assume that $G$ contains an unbalanced cycle. 
Note that by Lemma \ref{bnzl} it will suffice to show that $G$ has a balanced nowhere-zero $\mathbb{Z}_2 \times \mathbb{Z}_3$ flow.
To construct this, we fix an arbitrary orientation of $G$.   

Since $G$ has a nowhere-zero $\mathbb{Z}$-flow it cannot have a vertex of degree one.  Suppose that $G$ has a vertex $v$ of degree two.
If $v$ is incident with a loop, then this is the only edge and the result is trivial.  Otherwise, let $\delta(v) = \{e,f\}$, and note that by possibly replacing
$\sigma_G$ with an equivalent signature, we may assume $\sigma_G(e) = 1$.  By Lemma \ref{lred} and induction, the graph obtained by contracting $e$ has a balanced nowhere-zero $\mathbb{Z}_2 \times \mathbb{Z}_3$-flow $\phi$.  
Now we may extend the domain of $\phi$ to $E(G)$ by setting $\phi(e) = \pm \phi(f)$ so that $\phi$ is a balanced nowhere-zero
$\mathbb{Z}_2 \times \mathbb{Z}_3$-flow of $G$.  

If $G$ contains a vertex $v$ with ${\mathit deg}(v) \ge 4$, then by Lemma \ref{lred}, we may uncontract a positive edge at $v$ so that
the resulting graph $G'$ has a nowhere-zero $\mathbb{Z}$-flow.  By induction we may choose a balanced nowhere-zero 
$\mathbb{Z}_2 \times \mathbb{Z}_3$-flow $\phi$ of $G'$.  Now $\phi|_{E(G)}$ is a balanced nowhere-zero $\mathbb{Z}_2 \times \mathbb{Z}_3$-flow of $G$.
Thus, we may assume that $G$ is cubic.

Next suppose that there is a subset $X \subseteq V(G)$ with $|X| > 1$ so that $G[X]$ is
balanced and $|\delta(X)| \le 3$.  By possibly adjusting $\sigma_G$ and $\tau_G$, we may assume that 
$\sigma_G(e) = 1$ for every $e \in E(G)$ and that every half edge $h$ contained in an edge of
$\delta(X)$ and incident with a vertex in $X$ is directed toward this vertex (i.e. it satisfies $\tau_G(h) = 1$).  
Let $G_x$ be the graph obtained from $G$ by identifying $X$ to a single new
vertex $x$ and deleting any loops formed in this process.  By Lemma \ref{lred} and induction, we may choose a balanced nowhere-zero
$\mathbb{Z}_2 \times \mathbb{Z}_3$-flow $\phi^x$ of $G_x$.  Now starting with the bidirected graph $G$ 
we identify $V(G) \setminus X$ to a single new vertex $y$ (again deleting all newly formed loops) and then modify to get an ordinary directed graph 
$G_y$ by defining each edge in $\delta(y)$ to be positive and directed away from $y$.  
By Lemma \ref{L31} we may choose a nowhere-zero $\mathbb{Z}_2 \times \mathbb{Z}_3$-flow $\phi^y$ of $G_y$ so that
$\phi^y(e) = \phi^x(e)$ for every edge $e \in \delta(y)$.   Now the function 
$\phi : E(G) \rightarrow \mathbb{Z}_2 \times \mathbb{Z}_3$ given by the rule
\[ \phi(e) = \left\{ \begin{array}{ll}
    \phi^x(e) &   \mbox{if $e \in E(G_x)$}    \\
    \phi^y(e) &   \mbox{otherwise}
    \end{array} \right. \]
is a balanced nowhere-zero $\mathbb{Z}_2 \times \mathbb{Z}_3$-flow of $G$ as desired.  Thus, we may assume no such 
subset $X$ exists.  

If there is a balanced 4-cycle $C \subseteq G$, then we may assume $\sigma_G(e) = 1$ for every $e \in E(G)$.  
Let $G'$ be the graph obtained from $G$ by deleting $E(C)$ and then identifying $V(C)$ to a single new
vertex $v$.  By Lemma \ref{lred}, $G'$ has a nowhere-zero $\mathbb{Z}$-flow, so by induction we may choose a balanced nowhere-zero
$\mathbb{Z}_2 \times \mathbb{Z}_3$-flow $\phi$ of $G'$.  It is now straightforward to verify that $\phi$ can be 
extended to a balanced nowhere-zero $\mathbb{Z}_2 \times \mathbb{Z}_3$-flow of $G$ (if no edges incident with $v$ are in the support
of $\phi_1$ then we may extend $\phi$ so that $E(C)$ is in the support of $\phi_1$; otherwise we may 
choose $\phi$ so that at least two edges of $E(C)$ are in the support of $\phi_1$).

Since $G$ is not balanced, $G \setminus e$ must have an unbalanced cycle for every $e \in E(G)$, and it follows that $G$ is a shrubbery.  
If $G$ is 2-connected, then it must contain an unbalanced theta.  Otherwise, consider a leaf block of $G$.  Since this leaf block must contain an unbalanced cycle, it must either be an unbalanced loop edge, or it must contain an unbalanced theta.  In either case, we may apply 
Lemma \ref{workhorse} with $\epsilon = 1$ to choose a nowhere-zero watering of $G$.  This is precisely a balanced nowhere-zero 
$\mathbb{Z}_2 \times \mathbb{Z}_3$ flow, as desired.
\quad\quad$\Box$

\section{Removable Cycles}

In the proof of Seymour's 6-Flow Theorem, a key concept was that of a removable cycle.  In 
this section we introduce an analogous concept for shrubberies and show that it has the 
desired properties.  This is all in preparation for the proof of Lemma \ref{workhorse}.

Let $G$ be a signed graph.  We define $V_2(G) = \{ v \in V(G) | {\mathit deg}(v) = 2 \}$.  For a cycle
$C \subseteq G$ we let $\mathcal{U}(C)$ denote the set of chords $e$ of $C$ for which $C \cup e$ 
is an unbalanced theta, and we let $\mathcal{B}(C)$ denote those chords $e$ for which $C \cup e$ 
is a balanced theta.  We will call $C$ a {\it removable} cycle if it has one of the following properties.
\begin{enumerate}
\item $C$ is unbalanced.
\item $|V(C) \cap V_2(G)| + | {\mathcal U}(C) | \ge 2$.
\end{enumerate}

The following Lemma shows that removable cycles behave appropriately.  

\begin{lemma}
\label{L72}
Let $G$ be a shrubbery and let $C \subseteq G$ be a removable cycle.  Then,
for every nowhere-zero watering $\phi'$ of $G' = G \setminus V(C)$, there exists a
nowhere-zero watering $\phi$ of $G$ so that $\phi(e) = \phi'(e)$ for every
$e \in E(G')$ and 
${\mathit supp}(\phi_1) = E(C) \cup {\mathit supp}(\phi'_1)$.
\end{lemma}

\noindent{\bf Proof:}
Our first step will be to extend $\phi'$ to a function $\phi : E(G) \rightarrow \mathbb{Z}_2 \times \mathbb{Z}_3$ as follows.
Here we view $\alpha_e$ as a variable in $\mathbb{Z}_3$ for every $e \in \mathcal{U}(C)$.  
\[ \phi'(e) = \left\{ \begin{array}{ll}
    (0,\pm 1)&	\mbox{if $e \in \delta(V(C))$} \\
    (1,0)       &   \mbox{if $e \in E(C)$}          \\
    (0,1)       &   \mbox{if $e \in {\mathcal B}(C)$ }  \\
    (0,\alpha_e)    &   \mbox{if $e \in {\mathcal U}(C)$}   \\
    \end{array} \right. \]
Since every $v \in V(G) \setminus V(C)$ adjacent to a vertex in $V(C)$ has degree less than three in $G'$, we may 
choose values $\phi(e)$ for the edges $e \in \delta(V(C))$ so that $\phi$ satisfies the boundary condition for a watering
at every vertex in $V(G) \setminus V(C)$.  
By construction $\partial \phi_1(v) = 0$ for every $v \in V(C)$.  So we need only adjust $\partial \phi_2(v)$ for $v \in V(C)$ 
to obtain a watering.  We now split into cases based on $C$.

\bigskip

\noindent{\sl Case 1:} $C$ is unbalanced.
\smallskip

Choose arbitrary $\pm 1$ assignments to the variables $\alpha_e$.  Since $C$ is 
unbalanced, for every vertex $u \in V(C)$ there is a function $\eta^u : E(C) \rightarrow \mathbb{Z}_3$
so that $\partial \eta^u(v) = 0$ for every $v \in V(C) \setminus \{u\}$ and $\partial \eta^u(u) = 1$.  
Now we may adjust $\phi_2$ by adding a suitable combination of the $\eta^u$ functions so that $\phi$ 
is a nowhere-zero watering of $H$, as desired.

\bigskip

\noindent{\sl Case 2:} $C$ is balanced.

We may assume without loss that every edge of $C$ is positive and every unbalanced chord is oriented
so that each half edge is directed away from its end.  In this case each unbalanced chord $e$ contributes 
$-2 \phi_2(e) = \alpha_e$ to the sum $\sum_{v \in V(C)} \partial \phi_2 (v)$.  For every $v \in V(C) \cap V_2(G)$ 
let $\beta_v$ be a variable in $\mathbb{Z}_3$.  Since $| \mathcal{U}(C) | + |V(C) \cap V_2(G)| \ge 2$ we may now choose $\pm 1$ assignments to 
all of the variables $\alpha_e$ and $\beta_v$ so that the following equation is satisfied.
\[ \sum_{v \in V(C)} \partial \phi_2 (v) = \sum_{v \in V(C) \cap V_2(G)} \beta_v \]
Now by lemma \ref{posadj} we may choose a function $\psi : E(C) \rightarrow \mathbb{Z}_3$ so that 
\[ \partial \psi(v) = \left\{ \begin{array}{cl}
					\beta_v 	- \partial \phi_2(v) 		&	\mbox{if $v \in V(C) \cap V_2(G)$}	\\
							- \partial \phi_2(v) 		&	\mbox{if $v \in V(C) \setminus V_2(G)$}	
					\end{array} \right. \]	
Now modify $\phi$ by adding $\psi$ to $\phi_2$.  After this adjustment, $\phi$ is a nowhere-zero watering of $G$, as desired.
\quad\quad$\Box$

\section{Cycles in Signed Graphs}

The most difficult case for our approach to the proof of Lemma \ref{workhorse} is 
when $\epsilon = 1$ and $G$ is a cubic shrubbery without two disjoint odd cycles.  
In this case a watering $\phi$ satisfying the lemma must have the property that the support of 
$\phi_1$ consists only of balanced cycles.  So, we cannot take advantage of an unbalanced removable 
cycle.  Since there are no vertices of degree 2, our only hope for a suitable removable cycle is to find
a balanced cycle $C$ which has at least two unbalanced chords.  Through the course of this section
we will show that such a cycle does exist.  In fact we will prove a couple of more general results 
which will help in finding removable cycles in other cases too.

\begin{lemma}
\label{routeincubic}
Let $G$ be a 3-connected cubic graph, let $H \subset G$ be connected, and let $x,y \in V(H)$.  
Then there exists a path $P \subseteq H$ from $x$ to $y$ with the property that every component of $H \setminus E(P)$ 
contains a vertex $v$ with ${\mathit deg}_H(v) \le 2$
\end{lemma}

\noindent{\it Proof:} Choose a path $P \subseteq H$ from $x$ to $y$ according to the following criteria.  
\begin{enumerate}
\item Lexicographically maximize the sizes of the components of $H \setminus E(P)$ which contain a vertex of degree at most two in $H$.
\item Lexicographically maximize the sizes of the remaining components of $H \setminus E(P)$ (subject to condition 1.)
\end{enumerate}
So, our first priority is for the largest component of $H \setminus E(P)$ which contains a point in $V_2(H)$ to be as large as possible.  
Subject to this, our second priority is to maximize the size of the second largest component of $H \setminus E(P)$ which contains a 
vertex in $V_2(H)$ (if one exists) and so on.   

Suppose (for a contradiction) that $P$ does not satisfy the lemma and choose $H'$ to be the smallest component of 
$H \setminus E(P)$ which does not contain a vertex with degree at most two in $H$.  Let $P'$ be the minimal 
subpath of $P$ which contains all vertices in $V(H') \cap V(P)$.  If every vertex in $P'$ is in $H'$ then the ends of $P'$ form a 2-separation 
separating $H' \cup P'$ from the rest of the graph.  Since this is impossible, there must exist a vertex in the interior of $P'$ 
which is in another component of $H \setminus E(P)$.  Now consider rerouting the original path $P$ by replacing $P'$ with 
a path through $H'$ from one end of $P'$ to the other.  This new path contradicts our choice of $P$ (since some component 
other than $H'$ will get larger), giving us a contradiction.  
\quad\quad$\Box$

\bigskip

A subgraph $H \subseteq G$ is {\it peripheral} if $G \setminus V(H)$ is connected and no edge in $E(G) \setminus E(H)$ has both ends in $V(H)$.
Note that if $G$ is cubic and $H$ is a cycle, the above condition is equivalent to $G \setminus E(H)$ connected.  Tutte used the technique of lexicographic
maximization employed above to establish numerous important properties of peripheral cycles.  Here we will require the following of his results.

\begin{theorem}[Tutte \cite{Tu3}]
\label{tutte}
If $G$ is a 3-connected graph, the peripheral cycles of $G$ generate the cycle-space of $G$ over $\mathbb{Z}_2$. 
\end{theorem}


If $G$ is a signed cubic graph, we will say that a balanced cycle
$D \subseteq G$ is a {\it halo} if $G \setminus E(D)$ contains a pair
$(P_1,P_2)$ of vertex disjoint paths, called a {\it cross} of $D$, with
the following properties:

\smallskip
\begin{tabular}{lp{5.1in}}
(i)	& $P_i$ has both ends in $V(D)$ but is internally disjoint from $V(D)$ for $i=1,2$. \\ 
(ii)    &   $D \cup P_1 \cup P_2$ is isomorphic to a subdivision of $K_4$. \\
(iii)   &   $P_i \cup D$ contains an unbalanced cycle for $i=1,2$. \\
(iv)   &  Every component of $G \setminus E(D)$ contains either $P_1$ or $P_2$
\end{tabular}
\smallskip

\begin{lemma}
\label{L63}
Let $G$ be a 3-connected cubic signed graph which is not balanced.  If $G$ has a nowhere-zero $\mathbb{Z}$-flow
but does not have two vertex disjoint unbalanced cycles, then it contains a halo.
\end{lemma}

\noindent{\it Proof:}
It follows from Theorem \ref{tutte} that $G$ contains an unbalanced peripheral cycle $C$ (otherwise this theorem would force 
all cycles to be balanced).  By possibly replacing $\sigma_G$ with an equivalent signature, we may then assume that every edge in $E(G) \setminus E(C)$ 
is positive.  Let $E^- \subseteq E(C)$ be the set of negative edges.  Note that $|E^-|$ is odd (since $C$ is unbalanced) and $|E^-| > 1$ 
(since $G$ has  a nowhere-zero $\mathbb{Z}$-flow), so we must have $|E^-| \ge 3$.  

Choose a path $P \subseteq G \setminus E(C)$ so that the ends of $P$ are in distinct components of 
$C \setminus E^-$ and subject to this choose $P$ so as to lexicographically maximize the sizes of the components of the graph 
$G' = G \setminus (E(C) \cup E(P))$.  Note that since $C$ is peripheral, every component of $G'$ must contain a vertex in $P$.  
Let $A_0,B_0$ be the two components of $C \setminus E^-$ which contain an endpoint of $P$.  Call a component of $G'$ \emph{rich} if 
it contains a vertex in $V(C) \setminus (V(A_0) \cup V(B_0))$ and otherwise \emph{poor}

\begin{figure}[h]
  \centering
    \includegraphics[height=4.5cm]{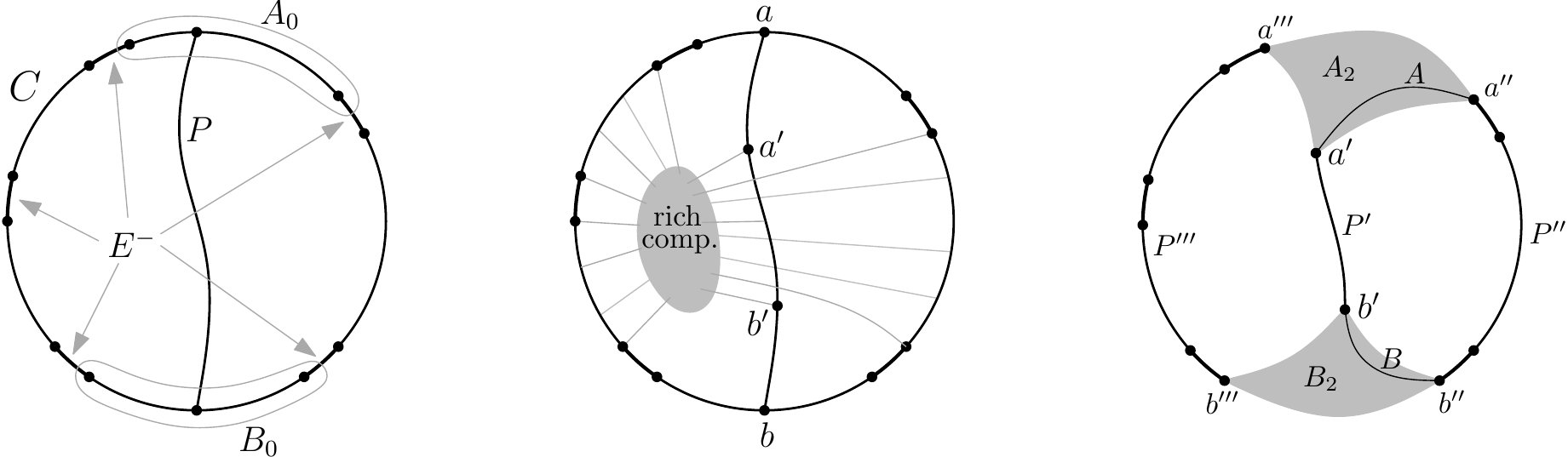}
  \label{halo_all}
\end{figure}

Suppose (for a contradiction) that there is a rich component $H$ of $G'$ and some other component $H'$ of equal or greater size.  Now 
$H$ contains a path from $P$ to a vertex in $V(C) \setminus (V(A_0) \cup V(B_0))$ so we may reroute $P$ using this path so as to 
increase the size of $H'$, thus contradicting our choice of $P$.  It follows that there is exactly one rich component of $G'$ and this is at least 
as large as any other component of $G'$.  

Let $a,b$ be the ends of $P$ in $A_0, B_0$.  Traversing the path $P$ from $a$ to $b$ we let $a'$ ($b'$) be the first (last) vertex encountered which is in the rich component of $G'$.  Define $A_1$ ($B_1)$ to be the union of $A_0$ ($B_0$) and the subpath of $P$ from $a$ to $a'$ ($b$ to $b'$).  If there exists a poor component of $G'$ which contains a vertex in $A_1$ and a vertex in $(C \cup P) \setminus A_1$ then by rerouting $P$ through this component we 
may increase the size of the rich component, thus contradicting the choice of $P$.   A similar argument for $B_1$ shows that every poor component
intersects $V(C \cup P)$ in either a subset of $V(A_1)$ or a subset of $V(B_1)$.  Define $A_2$ ($B_2$) to be the union of $A_1$ ($B_1$) and all 
poor components of $G'$ which intersect $V(C \cup P)$ in a subset of $V(A_1)$ ($V(B_1)$).  

Consider the unique path in $C$ from $a$ to $b$ which contains an even (odd) number of edges in $E^-$ and while traversing this path from $a$ 
to $b$ let $a''$ ($a'''$) be the last vertex in $A_0$ which is encountered, let $b''$ ($b'''$) be the first vertex of $B_0$ which is encountered, and let
$P''$ ($P'''$) be the subpath from $a''$ to $b''$ ($a'''$ to $b'''$).  Analogously define $P'$ to be the subpath of $P$ from $a'$ to $b'$.  Now apply lemma \ref{routeincubic} to choose a path $A$ in the subgraph $A_2$ from $a'$ to $a''$ and a path $B$ in $B_2$ from $b'$ to $b''$.  Define the cycle
$D = P' \cup A \cup P'' \cup B$.  It follows from our construction that one of the following is true.
\begin{enumerate}
\item $D$ is a peripheral cycle.
\item $P'''$ is a one edge path and $G \setminus E(D)$ consists of two components one of which is $P'''$.
\end{enumerate}
We claim that $D$ is a Halo.  As $D$ is balanced, to prove this we need only find suitable paths $P_1$ and $P_2$.  Since the 
path $P''$ contains a nonzero even number of edges in $E^-$ we may choose a vertex $u \in V(P'')$ so that the two maximal 
subpaths of $P''$ ending at $u$ both contain an odd number of edges in $E'$.  Now the unique rich component of $G'$ contains a 
path $P_1$ from $u$ to a vertex in $P'$.  Define $P_2$ to be the unique component of $C \setminus E(D)$ which contains all edges 
in $P'''$.  It follows from our construction that $(P_1,P_2)$ is a cross of $D$, as desired.
\quad\quad$\Box$

\bigskip

Let $D$ be a halo with cross $(P_1,P_2)$ and let $X$ be the set of vertices consisting of all endpoints of 
$P_1$ and $P_2$.  We define a \emph{side} to be a path in $D$ between two vertices in $X$ which is 
internally disjoint from $X$ (so there are four sides).  If $Q_1,Q_2$ are vertex disjoint sides, we call them 
\emph{opposite}.  The next lemma shows that under appropriate assumptions, we can modify a cross 
so that a pair of opposite sides both consist of a single edge, say $e_1$ and $e_2$.  At this point the graph consisting
of the halo together with the cross minus $e_1,e_2$ is a balanced cycle for which $e_1,e_2$ are unbalanced
chords.

\begin{lemma}
\label{L64}
Let $G$ be a signed 3-connected cubic graph without two disjoint unbalanced cycles.  
Let $D$ be a halo with cross $(P_1,P_2)$ and let $Q_1,Q_2$ be opposite sides.  Then there 
exists a cross $(P_1', P_2')$ of $D$ with opposite sides $Q_1',Q_2'$ satisfying
\begin{enumerate}
\item $|E(Q_i')| = 1$ for $i=1,2$.
\item $Q_i' \subseteq Q_i$ for $i=1,2$.
\end{enumerate}
\end{lemma}

\noindent{\bf Proof:}
Choose a cross $(P_1',P_2')$ of $D$ with opposite sides $Q_1',Q_2'$ of $C$ which satisfy the 
second property above, and subject to this $|E(Q_1')| + |E(Q_2')|$ is minimum.  If 
$|E(Q_1')| + |E(Q_2')| = 2$ then we are finished.  Otherwise, we may assume that $|E(Q_1')| \ge 2$, and we may
choose a vertex $v$ in the interior of $Q_1'$.  By the last property in the definition of Halo, we may choose a 
path $R \subseteq G \setminus E(D \cup P_1' \cup P_2')$ from $v$ to a vertex $u \in V(P_1' \cup P_2')$.  Without 
loss suppose that $u \in V(P_1')$ and let $w$ be the end of $P_1'$ contained in $Q_1'$.  Let $C$ be the 
cycle of $G$ consisting of the path in $P_1'$ from $u$ to $w$ the path in $Q_1'$ from $v$ to $w$ and the 
path $R$.  It follows from the assumption that $G$ does not have two disjoint unbalanced cycles that 
$C$ is balanced.  Now rerouting the path $P_1'$ using $R$ gives us a cross which contradicts the choice 
of $(P_1',P_2')$ thus completing the proof.
\quad\quad$\Box$

\section{Watering a Shrubbery}

In this section, we complete the proof of Lemma \ref{workhorse}.  This will require us to first establish 
one lemma on the existence of paths in shrubberies, which in turn will call upon the following
classical result (simplified to the case of subcubic graphs).

\begin{theorem}[Messner and Watkins \cite{MW}]
\label{meswat}
Let $G$ be a 2-connected graph with maximum degree three and let
$y_1,y_2,y_3 \in V(G)$.  Either there exists a cycle $C \subseteq G$ with $y_1,y_2,y_3 \in V(C)$ or there
is a partition of $V(G)$ into $\{X_1,X_2,Y_1,Y_2,Y_3\}$ with the following
properties:
\begin{enumerate}
\item $y_i \in Y_i$ for $1 \le i \le 3$.
\item There are no edges between $X_1$ and $X_2$ or between $Y_i$ and $Y_j$ for
$1 \le i < j \le 3$.
\item There is exactly one edge between $X_i$ and $Y_j$ for every $i=1,2$
and $j=1,2,3$.
\end{enumerate}
\end{theorem}

For a path $P$ we let ${\mathit Ends}(P)$ denote the ends of $P$ and ${\mathit Int}(P)$ the set of 
interior vertices of $P$.

\begin{lemma}
\label{L73}
Let $G$ be a balanced 2-connected shrubbery and let $v_1,v_2 \in V_2(G)$.
Then there is a path $P \subseteq G$ with ${\mathit Ends}(P) = \{v_1,v_2\}$ and 
$|{\mathit Int}(P) \cap V_2(G)| \ge 2$.
\end{lemma}

\noindent{\bf Proof:}
We proceed by induction on $|V(G)|$.  If there exists
$Y \subseteq V(G) \setminus \{v_1,v_2\}$ so that $G[Y]$ contains a cycle 
and $|\delta(Y)| = 2$, then choose a minimal set $Y$ with these
properties.  By construction, $G[Y]$ is 2-connected.  Let
$y_1,y_2 \in Y$ be the two vertices incident with an edge of $\delta(Y)$.
Inductively, we may choose a path $Q \subseteq G[Y]$ with
${\mathit Ends}(Q) = \{y_1,y_2\}$ and with
$|{\mathit Int}(Q) \cap V_2(G)| \ge 2$.  Since $G$ is 2-connected, we
may choose two vertex disjoint paths from $\{v_1,v_2\}$ to $\{y_1,y_2\}$ and these together with 
$Q$ give the desired path.  So, we may assume $G$ has no such set $Y$.

Choose vertices $y_1,y_2 \in V_2(G) \setminus \{v_1,v_2\}$ (these must exist by our definitions) and then modify $G$
to form a new graph $G'$ by adding a new vertex $y_3$ and joining $y_3$ to the
vertices $v_1,v_2$.  If there is a cycle $C \subseteq G'$ with $y_1,y_2,y_3 \in V(C)$, then $C \setminus y_3$ is a path of
$G$ which satisfies the lemma.  Otherwise Theorem \ref{meswat} gives us a partition
$\{ X_1,X_2,Y_1,Y_2,Y_3 \}$ of $V(G')$ with the given properties.  It follows from the above argument
that $G[Y_i]$ is a path for $i=1,2$.  If there is a vertex $w$ in 
$V_2(G) \setminus (\{y_1,y_2,v_1,v_2\} \cup Y_3 )$ then 
choose a cycle $C$ of $G'$ with $y_3,w \in V(C)$.  Now $C$ must contain either $y_1$ or $y_2$ 
so again $C \setminus y_3$ satisfies the lemma.  So, we may assume no such vertex exists.  
In this case, since $G$ is a shrubbery, we must have $|X_1| = |X_2| = |Y_1| = |Y_2| = 1$.
But this contradicts the assumption that $G$ has no balanced 4-cycle.
\quad\quad$\Box$

\bigskip

In a connected graph $G$ an edge cut $\delta(X)$ \emph{separates cycles} if both the subgraph induced on 
$X$ and the subgraph induced on $V(G) \setminus X$ contain cycles.  We say that $G$ is \emph{cyclically} 
$k$-\emph{edge connected} if every edge cut separating cycles has size at least $k$.  

\bigskip

\noindent{\bf Proof of Lemma \ref{workhorse}}
Suppose (for a contradiction) that the lemma is false, and choose a counterexample $G$ 
for which $|E(G)|$ is minimum.  We define $V_2 = V_2(G)$ and now proceed toward a contradiction with
a series of numbered claims.  

\bigskip
\noindent{(1) $G$ is 2-connected}
\smallskip

By the minimality of our counterexample, it follows immediately that $G$ is connected.  
Suppose (for a contradiction) that $G$ has a cut-edge $f$ and fix an arbitrary orientation of $G$.  
Let $G', G''$ be the components of $G \setminus f$ and note that by the minimality of our counterexample
we may choose nowhere-zero waterings $\phi',\phi''$ of $G', G''$.  If $G$ contains an unbalanced theta or a loop, then so
does $G'$ or $G''$, so in this case, we may also choose $\phi',\phi''$ so that
$\sigma_G( {\mathit supp}(\phi'_1) ) \sigma_G( {\mathit supp}(\phi''_1) ) = \epsilon$.
Next, choose $\alpha,\beta = \pm 1$ so that the function
$\phi : E(G) \rightarrow \mathbb{Z}_2 \times \mathbb{Z}_3$ given by
\[ \phi(e) = \left\{ \begin{array}{ll}
    \alpha \phi'(e)      &   \mbox{if $e \in E(G')$}  \\
    \beta  \phi''(e)      &   \mbox{if $e \in E(G'')$}  \\
    (0,1)           &   \mbox{if $e = f$}
    \end{array} \right. \]
is a nowhere-zero watering of $G$.  If $G$ contained an unbalanced theta or a loop, then $\sigma_G ( {\mathit supp} (\phi_1) ) = \epsilon$.  
This contradiction shows that (1) holds.

\bigskip
\noindent{(2) $G$ contains an unbalanced theta}
\smallskip

Suppose (for a contradiction) that (2) is false.  The result holds trivially when $G$ consists of a single loop edge, 
so by (1) $G$ has no loop.  In this case we only need to find a nowhere-zero watering of $G$.  
If $G$ contains an unbalanced cycle $C$, then this cycle is removable.  Otherwise $G$ is balanced, so $|V_2| \ge 4$ 
and again we have a removable cycle $C$.  By minimality $G \setminus V(C)$ has a nowhere-zero watering, and 
now we may use Lemma \ref{L72} to extend this to a nowhere-zero watering of $G$ (which is a contradiction).

\bigskip
\noindent{(3) $G$ does not contain a removable cycle with one of the following
properties}

(A) $G \setminus V(C)$ contains an unbalanced theta

(B) $G \setminus V(C)$ is  balanced, and $\sigma_G(C) = \epsilon$

\smallskip

Suppose (3) is violated and choose a nowhere-zero watering $\phi'$ of $G \setminus V(C)$.  In case (A), we may choose $\phi'$ so that
$\sigma_G ( {\mathit supp} (\phi'_1) ) = \epsilon \sigma_G(E(C))$.  Now by Lemma
\ref{L72}, we may extend $\phi'$ to a nowhere-zero watering $\phi$ of $G$ so
that ${\mathit supp}(\phi_1) = {\mathit supp}(\phi'_1) \cup E(C)$.  By
construction we have that $\sigma_G( {\mathit supp}(\phi_1) ) = \epsilon$.

\bigskip
\noindent{(4) $G$ does not contain two unbalanced cycles $C_1,C_2$ so that
$V(C_1) \cup V(C_2)$ contains all the vertices of degree three in $G$}
\smallskip

Again we suppose (4) fails and proceed toward a contradiction.  
If $G \setminus C_i$ contains an unbalanced theta, then $C_i$ is a removable
cycle satisfying (3A), so we are finished.  Thus, we may assume that every
component of $G \setminus (E(C_1) \cup E(C_2))$ is a path with one end in
$V(C_1)$ and the other end in $V(C_2)$.  If $\epsilon = -1$, then we may
choose an unbalanced cycle $C \subseteq G$ so that $G \setminus V(C)$ is a
forest.  This contradicts (3B).  So we may assume $\epsilon = 1$.  Now 
$G' = G \setminus V(C_1 \cup C_2)$ is a forest and we may choose a nowhere-zero
watering $\phi'$ of $G'$ with ${\mathit supp}(\phi_1) = \emptyset$.  By two applications of Lemma
\ref{L72}, we may extend $\phi'$ to a nowhere-zero watering $\phi$ of $G$ with
${\mathit supp}(\phi_1) = E(C_1) \cup E(C_2)$.  Now
$\sigma_G({\mathit supp}(\phi_1)) = 1$ as required.

\bigskip
\noindent{(5) There does not exist $X \subseteq V(G)$ so that $\delta(X)$ separates
cycles, $|\delta(X)| = 2$, and $G[ V(G) \setminus X ]$ contains an unbalanced theta.}
\smallskip

Suppose (5) is false and choose a minimal set $X$ with the above properties.  Since
$G[ V(G) \setminus X]$ contains an unbalanced theta, every removable cycle which is 
contained in $G[X]$ satisfies (3A), so it will suffice to show that such a cycle exists.  
If $G[X]$ contains an unbalanced cycle $C$, then $C$ is a removable cycle of $G$.  
Otherwise, we must have $|X \cap V_2| \ge 2$.  By the
minimality of $X$, the graph $G[X]$ is 2-connected.  Thus, we may choose a cycle
$C \subseteq G[X]$ with $|V(C) \cap V_2| \ge 2$.

\bigskip
\noindent{(6) There does not exist $X \subseteq V(G)$ so that $\delta(X)$ separates
cycles, $|\delta(X)| = 2$, and $G \setminus \delta(X)$ is balanced.}
\smallskip

Suppose (6) fails, choose a minimal set $X$ with the above properties, and let
$\delta(X) = \{e_1,e_2\}$.  By possibly replacing $\sigma_G$ by an equivalent
signature, we may assume that $\sigma_G(e_1) = -1$, and that $\sigma_G(e) = 1$
for every other edge $e \in E(G) \setminus \{e_1\}$ .  If
$\epsilon = -1$, then let $C$ be a cycle of $G$ with $e_1 \in E(C)$.  Then
$C$ is a removable cycle contradicting (3B) so we are done.  Thus, we may
assume that $\epsilon = 1$.  Now, $|X \cap V_2| \ge 2$ and by the
minimality of $X$, we have that $G[X]$ is 2-connected.  Thus, we may choose
a cycle $C \subseteq G[X]$ with $|V(C) \cap V_2| \ge 2$.
If $e_1$ is incident with a vertex of $V(C)$ or $e_1$ is a cut-edge of
$G \setminus V(C)$, then $C$ is a removable cycle satisfying (3B).  Otherwise,
$e_1$ is in an unbalanced theta of $G \setminus V(C)$, so $C$ is a removable
cycle satisfying (3A).

\bigskip
\noindent{(7) $G$ is cyclically 3-edge-connected}
\smallskip

Suppose (for a contradiction) that (7) is false and let $e \in E(G)$ be an edge in a 2-edge 
cut of $G$ which separates cycles.  Let $S = \{f \in E(G) \mid \mbox{ $\{e,f\}$ is an edge cut of $G$ }\} \cup \{e\}$,
and let $H_1,H_2,\ldots,H_m$ be the components of $G \setminus S$ which are not
isolated vertices.  
Note that $m \ge 2$.  By (5), we have that every $H_i$ is either 
balanced or it is an unbalanced cycle.  By (6) we may assume that $H_1$ is
an unbalanced cycle.  Let
$X_i = \{ v \in V(H_i) \mid \mbox{$v$ is incident with an edge in $S$} \}$ for
$1 \le i \le m$.  Now for every $2 \le i \le m$ we will choose a path
$P_i \subseteq H_i$ with ${\mathit Ends}(P_i) = X_i$ according to the following
strategy:  If $H_i$ is  balanced, then by Lemma \ref{L73} we may
choose  $P_i \subseteq H_i$ so that
$|V_2 \cap {\mathit Int}(P_i)| \ge 2$.  If $H_i$ is an unbalanced
cycle of size two, then let $P_i$ be a single edge path in $H_i$.  If $H_i$
is an unbalanced cycle of size at least three, then choose
$P_i \subseteq H_i$ so that
${\mathit Int}(P_i) \cap V_2 \neq \emptyset$.  Finally, choose
a path $P_1 \subseteq H_1$ so that ${\mathit Ends}(P_1) = X_1$ and so that
$C = \cup_{i=1}^m P_i \cup S$ is a cycle with $\sigma_G(C) = \epsilon$.
If one of $H_2,\ldots,H_m$ is  balanced, then
$|V_2 \cap V(C)| \ge 2$, so $C$ is removable.  Otherwise, by (4)
$m \ge 3$ so $|V_2 \cap V(C)| + |{\mathcal U}(C)| \ge 2$ and again
$C$ is removable.  In either case, $C$ contradicts (3B).

\bigskip
\noindent{(8) $G$ does not contain two disjoint unbalanced cycles}
\smallskip

If (8) is false for $C_1$ and $C_2$, then by (4) we may choose
a vertex $v \in V(G) \setminus (V(C_1) \cup V(C_2))$ of degree three.  By (7)
we may choose 3 internally disjoint paths from $v$ to $V(C_1) \cup V(C_2)$. 
We may assume without loss that at least two of these paths end in $V(C_2)$.  However, 
then $C_1$ is a removable cycle contradicting (3A). 

\bigskip
\noindent{(9) $\epsilon = 1$}
\smallskip

If $\epsilon = -1$, then by (2) we may choose an unbalanced cycle
$C \subseteq G$.  Now, $\sigma_G(C) = \epsilon$ and by (8) we have that
$G \setminus C$ is balanced.  Thus, $C$ is a removable cycle contradicting
(3B).

\bigskip
\noindent{(10)} There is no edge $e = y_1 y_2 \in E(G)$ with $G' = G \setminus e$ balanced 
and ${\mathit deg}(y_i) = 3$ for $i=1,2$.

\smallskip

Suppose (for a contradiction) that (10) is false.  
It follows from the second property in the definition of shrubbery that $G$ contains a vertex $y_3 \in V_2$.  
If $G'$ contains a cycle $C$ with $\{y_1,y_2,y_3\} \in V(C)$, then $C$ is a removable cycle of $G$ contradicting (3B).
Otherwise we may choose a partition of $V(G')$ into $\{ X_1, X_2, Y_1, Y_2, Y_3 \}$ in accordance with 
Theorem \ref{meswat}.  Note that by (7) $G[Y_3]$ is a path.  If there is a vertex
$w \in V_2 \setminus \{y_3\}$, then any cycle $C \subseteq G'$ with
$w,y_3 \in V(C)$ is a removable cycle of $G$ contradicting (3B).  Therefore 
$V_2 = \{y_3\}$.  Since $G[X_i]$ is balanced for $i=1,2$ and $G[Y_j]$ is balanced for $j=1,2,3$ 
it then follows that $|X_i| = 1$ for $i=1,2$ and $|Y_j| = 1$ for $j = 1,2,3$.  However, then $G$ contains 
a balanced cycle of length four, contradicting our assumption.  

\bigskip
\noindent{(11) There is no edge $e \in E(G)$ so that $G \setminus e$ is
 balanced}
\smallskip

Suppose (11) does not hold and let $P \subseteq G$ be a maximal path of $G$ for which 
$e \in E(P)$ and ${\mathit Int}(P) \subseteq V_2$.  By (10) we may assume
that $|E(P)| \ge 2$.  Let ${\mathit Ends}(P) = \{y_0, y_1\}$ and let
$G' = G \setminus {\mathit Int}(P)$.  Now, $G'$ is 
balanced, so we may choose
$\{y_2,y_3\} \subseteq V_2 \setminus \{y_0,y_1\}$.  If $G'$
has a cycle $C$ containing $y_1,y_2,y_3$ then this is a removable cycle of $G$ which 
contradicts (3B).  Otherwise, we may choose a partition $\{ X_1, X_2, Y_1, Y_2, Y_3 \}$ 
in accordance with Theorem \ref{meswat}.  If $y_0 \in Y_2 \cup Y_3$ then we may choose a cycle of $G'$
containing $y_2$ and $y_3$, and this contradicts (3B).  It then follows from the cyclic connectivity
of $G$ that both $G[Y_2]$ and $G[Y_3]$ are paths.  
If there is a vertex $w$ in $V_2 \setminus (\{y_0,y_2,y_3\} \cup Y_1)$ then we may choose a cycle of $G'$ containing
$w$ and $y_1$, and since this cycle will also contain one of $y_2$ or $y_3$ it contradicts (3B).  If $y_0 \in Y_1$ then
$X_1$ and $X_2$ have no vertices of degree 2 so we have $|X_1| = |X_2| = |Y_2| = |Y_3| = 1$ giving us a balanced 
4-cycle, which is contradictory.  So, we may assume (without loss) that $y_0 \in X_1$.
For $i=2,3$ let $x_i$ be the unique vertex in $X_1$ with a neighbour in $Y_i$.  If $x_2 \neq x_3$ 
and $G[X_1]$ has a path from $x_2$ to $x_3$ which has $y_0$ as an interior vertex, then 
$G'$ has a cycle containing $y_0,y_2,y_3$ and we have a contradiction to (3B).  Otherwise 
there is a partition of $X_1$ into $\{X_1', X_1''\}$ so that $y_0 \in X_1'$ and $x_2,x_3 \in X_1''$ and there is just 
one edge between $X_1'$ and $X_1''$.  However in this case our original graph $G$ has a two edge cut $\delta(X_1' \cup Y_1)$ which
separates cycles which contradicts (7).

\bigskip

Let $\widehat{G}$ be the graph obtained from $G$ by suppressing all vertices of degree two.  Now 
every edge $e \in E(\widehat{G})$ is associated with a path $P_e$ of $G$ with the property that
${\mathit Int}(P_e) \subseteq V_2$ and ${\mathit Ends}(P_e) \cap V_2 = \emptyset$.  Define a 
signature $\widehat{\sigma}$ of $\widehat{G}$ by setting $\widehat{\sigma}(e) = \sigma(P_e)$ for every $e \in \widehat{G}$.  
By (7) $\widehat{G}$ is a 3-connected cubic graph with signature $\widehat{\sigma}$, and by (11), $\widehat{G} \setminus e$ contains an
unbalanced cycle for every edge $e$.  Thus, by Lemma \ref{L63} we may choose a halo $\widehat{C} \subseteq \widehat{G}$.  
Fix a cross of $\widehat{C}$ and let $\widehat{A}_1,\widehat{A}_2,\widehat{A}_3,\widehat{A}_4$ be the sides of $\widehat{C}$ with respect 
to this cross with $\widehat{A}_1$ and $\widehat{A}_3$ opposite.  Let $A_1,A_2,A_3,A_4$ be the corresponding paths of $G$ and assume that
$|V_2 \cap V(A_1 \cup A_3)| \le |V_2 \cap V(A_2 \cup A_4)|$.  
Now, by Lemma \ref{L64} we may choose a cross $(\widehat{P}_1, \widehat{P}_2)$ of $\widehat{C}$ with sides 
$\widehat{B}_i$ for $1 \le i \le 4$ so that  $\widehat{B}_j \subseteq \widehat{A}_j$ and $|E(\widehat{B}_j)| = 1$ for $j=1,3$.  
Let $P_i$ for $i=1,2$ and $B_j$ for $1 \le j \le 4$ 
be the paths of $G$ corresponding to $\widehat{P}_i$ and $\widehat{B}_j$.  We claim that the cycle 
$D = P_1 \cup P_2 \cup B_2 \cup B_4$ is removable.  To see this, observe that for $j = 1,3$ whenever 
$V_2 \cap V(B_j) = \emptyset$ the path $B_j$ consists of a single edge which is an unbalanced chord of $D$.  
Since $|V_2 \cap V(B_1 \cup B_3)| \le |V_2 \cap V(B_2 \cup B_4)| \le |V_2 \cap V(D)|$ we have 
$\mathcal{U}(D) + |V_2 \cap V(D)| \ge 2$ as required.  Therefore the cycle $D$ contradicts (3B) and 
this completes the proof.
\quad\quad$\Box$

\section*{Acknowledgement}

We would like to thank Paul Seymour for numerous fruitful discussions.

\end{document}